\documentclass[11pt, twoside]{amsproc}

\usepackage{amsmath}
\usepackage{amssymb, latexsym, mathrsfs}
\usepackage{graphicx, subfigure}
\usepackage{enumitem, endnotes}
\usepackage[usenames, dvipsnames]{color}
\usepackage[colorlinks=true,
	raiselinks=true,
	pagecolor=magenta,
	linkcolor=myblue,
	citecolor=mygreen,
	urlcolor=mygreen,
	pdfauthor=Debasish Chatterjee,
	pdftitle={On Stability of Randomly Switched Nonlinear Systems},
	pdfkeywords={switched systems, random switching, stochastic stability, stabilization},
	pdfsubject={Research Article},
	plainpages=false]{hyperref}

\definecolor{myred}{rgb}{0.6, 0, 0}
\definecolor{mygreen}{rgb}{0, 0.5, 0}
\definecolor{myblue}{rgb}{0, 0, 0.5}
\definecolor{mycyan}{rgb}{0, 0.5, 0.5}

\definecolor{myred}{rgb}{0.75, 0, 0}
\definecolor{mygreen}{rgb}{0, 0.75, 0}
\definecolor{myblue}{rgb}{0, 0, 0.75}
\definecolor{mycyan}{rgb}{0, 0.5, 0.5}

\newcommand{\R}{\ensuremath{\mathbb{R}}}
\newcommand{\N}{\ensuremath{\mathbb{N}}}
\newcommand{\PSet}{\ensuremath{\mathcal{P}}}

\newcommand{\ra}{\ensuremath{\rightarrow}}
\newcommand{\Ra}{\ensuremath{\;\Longrightarrow\;}}
\newcommand{\lra}{\ensuremath{\longrightarrow}}

\newcommand{\posR}{\ensuremath{\R_{\ge 0}}}
\newcommand{\norm}[1]{\ensuremath{\left\lVert #1 \right\rVert}}

\newcommand{\eps}{\ensuremath{\varepsilon}}
\newcommand{\eqnref}[1]{{\rm (\ref{#1})}}
\newcommand{\fa}{\ensuremath{\forall\,}}

\newcommand{\RemarkEnd}{\hspace{\stretch{1}}{$\vartriangleleft$}}

\newcommand{\DefEnd}{\hspace{\stretch{1}}{$\Diamond$}}

\newcommand{\ClassK}{\ensuremath{\mathcal{K}}}
\newcommand{\ClassKinfty}{\ensuremath{\mathcal{K}_{\infty}}}
\newcommand{\ClassKL}{\ensuremath{\mathcal{KL}}}

\renewcommand{\le}{\ensuremath{\leqslant}}
\renewcommand{\ge}{\ensuremath{\geqslant}}

\renewcommand{\limsup}{\ensuremath{\varlimsup}}

\newcommand{\epower}[1]{\ensuremath{\mathrm{e}^{#1}}}
\newcommand{\cardP}{\ensuremath{\mathrm{N}}}
\newcommand{\transp}{\ensuremath{^{\scriptscriptstyle{\mathrm T}}}}
\newcommand{\sigalg}{\ensuremath{\mathfrak{F}}}
\newcommand{\ProbSpace}[1]{\ensuremath{\mathbf{\Omega}}}
\newcommand{\Expec}[1]{\ensuremath{\mathsf{E}\!\left[#1\right]}}
\newcommand{\Prob}[1]{\ensuremath{\mathsf{P}\!\left(#1\vphantom{\big)}\right)}}

\newcommand{\xz}{\ensuremath{x_0}}

\newcommand{\drv}{\ensuremath{\mathrm{d}}}
\newcommand{\ol}{\overline}

\newcommand{\wt}{\widetilde}

\newcommand{\gasas}{{\sc gas}~a.s.}
\newcommand{\gasm}{{\sc gas-m}}

\newcommand{\idmat}[1]{\ensuremath{{I}_{#1\times#1}^{\vphantom{T}}}}

\newcommand{\therex}{\ensuremath{\exists\,}}
\newcommand{\setmin}{\ensuremath{\!\smallsetminus\!}}
\newcommand{\LieD}[2]{\ensuremath{\mathrm L_{#1}{#2}}}
\newcommand{\indic}[1]{\ensuremath{\boldsymbol{1}_{#1}}}

\numberwithin{equation}{section}
\swapnumbers
\newtheorem{theorem}[equation]{Theorem}
\newtheorem{corollary}[equation]{Corollary}
\newtheorem{lemma}[equation]{Lemma}

\theoremstyle{definition}
\newtheorem{defn}[equation]{Definition}
\theoremstyle{remark}
\newtheorem{remark}[equation]{Remark}

\theoremstyle{remark}

\title{On Stability of Randomly Switched Nonlinear Systems}
\author{Debasish Chatterjee}
\address{Physikstrasse 3, ETL K24, ETH Zurich, 8092 Zurich, Switzerland}
\email{chatterd@control.ee.ethz.ch}
\author{Daniel Liberzon}
\address{144 Coordinated Science Laboratory, University of Illinois at Urbana-Champaign, Urbana, IL 61801.}
\email{liberzon@uiuc.edu}

\begin{document}

    \begin{abstract}
		This article is concerned with stability analysis and stabilization of randomly switched nonlinear systems. These systems may be regarded as piecewise deterministic stochastic systems: the discrete switches are triggered by a stochastic process which is independent of the state of the system, and between two consecutive switching instants the dynamics are deterministic. Our results provide sufficient conditions for almost sure global asymptotic stability using Lyapunov-based methods when individual subsystems are stable and a certain ``slow switching'' condition holds. This slow switching condition takes the form of an asymptotic upper bound on the probability mass function of the number of switches that occur between the initial and current time instants. This condition is shown to hold for switching signals coming from the states of finite-dimensional continuous-time Markov chains; our results therefore hold for Markov jump systems in particular. For systems with control inputs we provide explicit control schemes for feedback stabilization using the universal formula for stabilization of nonlinear systems.
    \end{abstract}

	\keywords{switched systems, random switching, almost sure and mean stochastic stability, stabilization}

    \maketitle
	\section{Introduction}
		\label{s:intro}
		Randomly switched systems consist of a family of subsystems, together with a random switching signal which specifies the active subsystem at every instant of time. Since the dynamics are governed by an ordinary differential equation between any two successive switching instants, these systems may be regarded as piecewise deterministic stochastic systems~\cite{ref:davisMarkovModelsOptbk}. These systems have variable structure, and can be used as models for systems affected by random structural changes. Applications of randomly switched systems include economic and manufacturing systems, communication and biological systems affected by random delays and component failures, etc. One particularly interesting phenomenon is observed in certain sea snails ``Pleurobranchea'' and ``Tritonia.''\footnote{This particular phenomenon was communicated to the authors by P.~\!K.~\!Jha.} These organisms have simple neural networks and persistent stimuli cause them to swim. Random changes in stimuli, e.g., scent of random food locations or random noxious environmental conditions, cause them to take orienting turns towards food, or avoidance turns away from the noxious agents, respectively.

		A particular class of piecewise deterministic stochastic systems has received widespread attention, namely, Markov jump linear systems (MJLS). These systems may be realized as a family of linear subsystems, together with a switching signal generated by the state of a continuous-time Markov chain. Stability and stabilization (see~\cite{ref:kozinsurvey} for a detailed survey on different notions and results on stochastic stability) of MJLS have been extensively investigated, especially under the assumption that the parameters of the Markov chain are completely known, see e.g.~\cite{ref:bolzernCDC04, ref:jichizeck90, ref:fengMarkovJump92, ref:maoexpstabstocdelay} and the references therein. In particular, almost sure stabilization and stabilization in the mean of MJLS is discussed in~\cite{ref:fengMarkovJump92}, where the authors also establish precise equivalences between different stability notions for MJLS.

		In this article we neither restrict ourselves to linear subsystems nor to Markov switching signals. Our results provide sufficient conditions for almost sure stability of randomly switched nonlinear systems when each subsystem is stable, and the switching is ``slow'' in a certain statistical sense. The slow switching condition takes the form of an upper bound on the probability mass function of the number of switches between the initial and current time instants. This condition is shown to hold in the case of switching signals coming from finite-dimensional continuous-time Markov chains; consequently our results can be applied to Markov jump systems under appropriate conditions on the parameters of the generator matrices of the underlying Markov chains. Since almost sure stability implies stability in probability~\cite{ref:khasminskiiStochastic}, our results also provide sufficient conditions for stability in probability of randomly switched systems; a variant of stability in the mean is also obtained. Based on our analysis, we propose control schemes which achieve almost sure stabilization and stabilization in the mean for systems with control inputs, by employing the universal formula for nonlinear feedback stabilization~\cite{ref:sontagunivformula}.

		A myriad of techniques have been employed to study stability and stabilization of piecewise deterministic stochastic systems. HJB-based optimal control schemes for piecewise deterministic stochastic systems are well-studied, see e.g.~\cite{ref:davisMarkovModelsOptbk} for a detailed account. Linear control systems admit analytically solvable linear quadratic optimal design methods, and such techniques have been effectively combined with the stochastic nature of structural variations in~\cite{ref:jichizeck90}; stabilization schemes based on Lyapunov exponents are studied in~\cite{ref:fengMarkovJump92}. Game-theoretic techniques~\cite{ref:basarMinMaxSw} in the presence of disturbance inputs, and spectral theory of Markov operators~\cite{ref:meynODESpecMarkovOp} have also been employed to study these systems. Stabilization schemes using robust control methods are investigated in~\cite{ref:xiong05}; see also the references cited there. Stochastic hybrid systems, where the switching signal and its transition probabilities are state-dependent, are studied in~\cite{ref:lygerosGenStochHybSys, ref:hespanhaCommNet}, using an extended definition of the infinitesimal generator and optimal control strategies, respectively.

		In contrast to the above references, our techniques rely on the method of multiple Lyapunov functions, discussed in the context of stability analysis and stabilization of deterministic switched systems in~\cite[Chapter~3]{ref:liberzonbk}. Our results conceptually parallel the ones on deterministic switched systems, which require stability of individual subsystems and a slow switching condition, see e.g.~\cite[Chapter~3]{ref:liberzonbk}. We employ multiple Lyapunov functions for stability analysis and stabilizing controller design, coupled with suitable assumptions to take care of the stochastic nature of the switching signal. Recently a method of stabilization in probability of Markov jump systems, with control and white noise disturbance inputs for each subsystem, has been proposed in~\cite{ref:battilottiDwellTimeMJS}, which is similar in spirit. We propose stronger results that apply to a wider class of systems and switching signals, although our model is simpler due to the absence of noise.

    \section{Preliminaries}
        \label{s:Prelims}
        Let the Euclidean norm be denoted by $\norm\cdot$, the interval $[0, \infty[$ by $\posR$, and the set of natural numbers $\{1, 2, \ldots\}$ by $\N$. Recall that a continuous function $\alpha:\posR\lra\posR$ is of class $\ClassK$ if $\alpha$ is strictly increasing with $\alpha(0) = 0$, of class $\ClassKinfty$ if in addition $\alpha(r) \ra \infty$ as $r \ra \infty$; we write $\alpha\in\ClassK$ and $\alpha\in\ClassKinfty$ respectively.

        We define the family of systems
        \begin{equation}
            \dot x = f_p(x),\qquad p\in\PSet,
            \label{e:ssysfam}
        \end{equation}
        where the state $x\in\R^n$, $\PSet$ is a finite index set of $\cardP$ elements: $\PSet = \{1, \ldots, \cardP\}$, the functions $f_p:\R^n\lra\R^n$ are locally Lipschitz, $f_p(0) = 0$, $p\in\PSet$. 

		Let $(\Omega, \sigalg, (\sigalg_t)_{t\ge 0}, {\mathsf P})$ be a complete filtered probability space, where $(\sigalg_t)_{t\ge 0}$ be the natural filtration (satisfying the usual conditions) generated by a c\'adl\'ag random process $\sigma$ taking values in $\PSet$. We let $\sigma$ generate the \emph{randomly switched system}
        \begin{equation}
            \dot x = f_\sigma(x), \qquad x(0) = \xz, \quad t\ge 0.
            \label{e:ssysdef}
        \end{equation}
		from the family~\eqref{e:ssysfam}. We assume that there are no jumps in the state $x$ at the switching instants. Let the switching instants of $\sigma$ be denoted by $\tau_i, \;i = 1, 2, \ldots$, and $\tau_0 := 0$ by convention. As a consequence of the hypotheses of our main result, there is no explosion almost surely (see Remark~\ref{r:noexpl} for details); therefore the sequence $(\tau_i)_{i\ge 0}$ is divergent. Finally, we assume that for every compact subset $K\subset\posR\times\R^n$ there exists an integrable function $m_K^{\vphantom{K}}$ satisfying $\sup_{p\in\PSet}\norm{f_p(x)} \le m_K^{\vphantom{K}}(t)$ for all $(t, x)\in K$. Hence almost surely there exists a unique solution to~\eqref{e:ssysdef} in the sense of Carath{\'e}odory~\cite{ref:filippovbk} over a nontrivial time interval containing $0$; existence and uniqueness of a global solution will follow from the hypotheses of our main result. We let $x(\cdot)$ denote this solution. For $\xz = 0$, the solution to~\eqref{e:ssysdef} is identically $0$ for every $\sigma$; we shall ignore this trivial case in the sequel.

		We focus on the notion of stability defined below.
		\begin{defn}
			The system~\eqref{e:ssysdef} is said to be \textbf{\textit{globally asymptotically stable almost surely}} (\gasas{}) if the following two properties are simultaneously verified:
			\begin{enumerate}[label={\rm (SP\arabic*)}, align=right, leftmargin=*]
				\item $\mathsf P\!\left(\vphantom{\displaystyle{\sum}}\fa\eps > 0\;\; \therex\delta > 0\right.$ $\norm{\xz} < \delta\Ra\sup_{t\ge 0}\norm{x(t)} < \eps\left.\vphantom{\displaystyle{\sum}}\right) = 1$;
				\item $\mathsf P\!\left(\vphantom{\displaystyle{\sum}}\fa r, \eps' > 0\;\; \therex T \ge 0\right.$ $\norm{\xz} < r\Ra\sup_{t\ge T}\norm{x(t)} < \eps'\left.\vphantom{\displaystyle{\sum}}\right) = 1$.\DefEnd
			\end{enumerate}
			\label{d:gasas}
		\end{defn}

		This definition can be readily recast into the $\ClassKL$ framework as follows. The system~\eqref{e:ssysdef} is {\sc gas}~a.s.\ if
		\[
			\Prob{\therex\beta\in\ClassKL\;\;\fa\xz\in\R^n\;\;\fa t\ge 0 \quad\norm{x(t)}\le \beta(\norm{\xz}, t)} = 1.
		\]
		Notice that there is no assertion about uniformity of trajectories; it merely states that the state trajectory corresponding to almost every event satisfies some class-$\ClassKL$ bound.

	\section{Stability under Random Switching}
	\label{s:mainres}
	\subsection{Global asymptotic stability almost surely}
		For a switching signal $\sigma$, we denote the number of switches on the interval $]0, t]$ by $N_\sigma(t)$. The following main result of this article provides sufficient conditions for \gasas{} of~\eqref{e:ssysdef}.
		\begin{theorem}
			Consider the system~\eqref{e:ssysdef}. Suppose that there exist continuously differentiable functions $V_p:\R^n\lra\posR,\;p\in\PSet$, functions $\alpha_1, \alpha_2\in\ClassKinfty$, and real numbers $\wt\lambda, \ol\lambda, \lambda_\circ > 0$, $\mu > 1$, such that
			\begin{enumerate}[label={\rm (\roman*)}, align=right, leftmargin=*]
				\item $\alpha_1(\norm x) \le V_p(x) \le \alpha_2(\norm x)\quad \fa x\in\R^n,\;\fa p\in\PSet$;\label{cond:randstab:1}\medskip
				\item $\displaystyle{\frac{\partial V_p}{\partial x}f_p(x) \le -\lambda_\circ V_p(x)}\quad\fa x\in\R^n,\;\fa p\in\PSet$;\label{cond:randstab:2}\medskip
				\item $V_{p_1}(x) \le \mu V_{p_2}(x)\quad\fa x\in\R^n,\;\fa p_1, p_2\in\PSet$;\label{cond:randstab:3}
				\item $\therex M\in\N\cup\{0\}$ such that $\fa k\ge M$ we have $\displaystyle{\Prob{N_\sigma(t) = k} \le \frac{(\ol\lambda t)^k\epower{-\wt\lambda t}}{k!}}$;\label{cond:randstab:4}
				\item $\mu < \left(\lambda_\circ+\wt\lambda\right)/\ol\lambda$.\label{cond:randstab:5}
			\end{enumerate}
			Then~\eqref{e:ssysdef} is \gasas{}
			\label{t:randstab}
		\end{theorem}

		Before proving Theorem~\ref{t:randstab}, let us make the following observations.

		\begin{remark}
			Hypothesis~\eqref{cond:randstab:2} of Theorem~\ref{t:randstab} implies that each subsystem of the family~\eqref{e:ssysfam} is globally asymptotically stable; the right-hand side of the inequality being linear in $V_p$ is no loss of generality, see~\cite[Theorem~2.6.10]{ref:lakshmikanthamDifferentialIntegralInequalities1}. Hypothesis~\eqref{cond:randstab:3} first appeared in~\cite{ref:hespanhaAvgDwellTime}, and has almost become a standard in deterministic switched systems literature. It restricts the class of applicable Lyapunov functions by requiring the existence of a maximal global constant ratio among the functions, but it is not known whether this hypothesis actually incurs a loss of generality. Quadratic Lyapunov functions are universally considered for linear systems, and in this case the existence of a global constant $\mu$ is automatically guaranteed. Since the left hand side of the inequality in hypothesis~\ref{cond:randstab:4} is a probability measure of an event, the right hand side may be replaced by $\bigl(\epower{-\wt\lambda t}(\ol\lambda t)^k/k!\bigr)\wedge 1$. In the special case when $\wt\lambda = \ol\lambda$, this hypothesis states that the number of switches $N_\sigma(t)$ on $]0, t]$ of $\sigma$ is eventually upper bounded by the probability mass function of a Poisson process of parameter $\ol\lambda$.\RemarkEnd
			\label{r:lambdarelations}
		\end{remark}

		\begin{remark}
			Suppose $\sigma$ satisfies hypothesis~\ref{cond:randstab:4} of Theorem~\ref{t:randstab}. Then the probability of an explosion, i.e., of an accumulation of infinitely many switches over a finite time interval is zero. Indeed, if $\zeta\in\posR$, the event that there is an explosion at $t = \zeta$ is given by $\bigcap_{\nu\in\N}\big\{N_\sigma(\zeta) \ge \nu\big\}$. But $\Prob{\bigcap_{\nu\in\N}\big\{N_\sigma(\zeta) \ge \nu\big\}} \le \limsup_{\nu\uparrow\infty}\sum_{k = \nu}^\infty\Prob{N_\sigma(\zeta) = k}$, and using hypothesis~\ref{cond:randstab:4}, we get $\limsup_{\nu\uparrow\infty}\sum_{k = \nu}^\infty\Prob{N_\sigma(\zeta) = k} \le \limsup_{\nu\uparrow\infty}\sum_{k = \nu}^\infty\epower{-\wt\lambda\zeta}\frac{(\ol\lambda\zeta)^k}{k!} = 0$. Since $\zeta\in\posR$ is arbitrary, we conclude that there is no explosion almost surely.\RemarkEnd
			\label{r:noexpl}
		\end{remark}

The proof of Theorem~\ref{t:randstab} follows the sequence of Lemmas below. The idea of the proof may be briefly summarized as follows. The property (SP2) is proved first by first estimating the expected value of $V_{\sigma(t)}(x(t))$ for an arbitrary $t\ge 0$ via the moment generating function of $N_\sigma(t)$, and then proving a.s. asymptotic convergence of $\bigl(V_{\sigma(t)}(x(t))\bigr)_{t\ge 0}$ via Tonelli's theorem and an auxiliary Lemma that proves asymptotic convergence of $\norm{x(\cdot)}$ from the finiteness of a certain nonnegative integral. We also observe that since the (finite) family of subsystems is uniformly locally Lipschitz, the maximal temporal growth rate of trajectories is upper bounded by a constant in a neighborhood of $0$. The (SP1) property can now be established utilizing this fact and the (SP2) property that is proved first, thereby completing the proof.

		\begin{lemma}
			Suppose that hypotheses~\ref{cond:randstab:2} and~\ref{cond:randstab:3} of {\rm Theorem~\ref{t:randstab}} hold. Then we have $\Expec{V_{\sigma(t)}(x(t))} \le \Expec{\epower{(\ln\mu)N_\sigma(t)}}V_{\sigma(0)}(\xz)\epower{-\lambda_\circ t}\quad \fa t\ge 0$.
			\label{l:Vsigmaest}
		\end{lemma}
		\begin{proof}
			Recall that $(\tau_i)_{i\in\N}$ are the switching instants of $\sigma$. It follows from hypothesis~\ref{cond:randstab:2} that for $t\in[\tau_i, \tau_{i+1}[$, we have
			\[
				V_{\sigma(\tau_i)}(x(t)) \le V_{\sigma(\tau_i)}(x(\tau_i))\epower{-\lambda_\circ(t-\tau_i)}.
			\]
			In conjunction with hypothesis~\ref{cond:randstab:3}, this yields
			\[
				V_{\sigma(\tau_{i+1})}(x(\tau_{i+1})) \le \mu V_{\sigma(\tau_i)}(x(\tau_i))\epower{-\lambda_\circ(\tau_{i+1}-\tau_i)}.
			\]
			Iterating the last inequality from $i = 0$ to $i = N_\sigma(t)$ for an arbitrary time $t > 0$, we arrive at
			\begin{equation*}
				V_{\sigma(t)}(x(t)) \le \mu^{N_\sigma(t)}\epower{-\lambda_\circ t}V_{\sigma(0)}(\xz).
			\end{equation*}
			Since the initial condition is deterministic, taking expectations on both sides of the above inequality we get
			\begin{equation*}
				\Expec{V_{\sigma(t)}(x(t))} \le \Expec{\mu^{N_\sigma(t)}}\epower{-\lambda_\circ t}V_{\sigma(0)}(\xz),
			\end{equation*}
			which proves the claim.
		\end{proof}

		\begin{lemma}
			Suppose that hypothesis~\ref{cond:randstab:4} of {\rm Theorem~\ref{t:randstab}} holds. Then $\therex S \ge 0$ such that the {\em moment generating function} $\Expec{\epower{sN_\sigma(t)}}$ of $N_\sigma(t)$ satisfies $\Expec{\epower{sN_\sigma(t)}} \le S + \epower{\left(\epower{s}\ol\lambda - \wt\lambda\right)t}\quad \fa s\ge 0$.
			\label{l:mgenfun}
		\end{lemma}
		\begin{proof}
			Using hypothesis~\ref{cond:randstab:4}, for $s\ge 0$ a little computation leads to
			\begin{align*}
				& \Expec{\epower{sN_\sigma(t)}} = \sum_{k = 0}^\infty \epower{sk}\Prob{N_\sigma(t) = k}\nonumber\\
				& \le \sum_{k = 0}^{M-1} \epower{sk}\Prob{N_\sigma(t) = k} + \sum_{k = M}^\infty \epower{sk}\frac{(\ol\lambda t)^k\epower{-\wt\lambda t}}{k!}\nonumber\\
				& \le S + \epower{\left(\epower{s}\ol\lambda-\wt\lambda\right)t},
			\end{align*}
			where $S := \sum_{k = 0}^{M-1}\epower{sk} \ge 0$. Clearly, $\Expec{\epower{sN_\sigma(t)}}$ is well defined for $t\ge 0$. 
		\end{proof}

		The following Barbalat-type lemma will be needed to conclude asymptotic convergence of the state trajectory from the finiteness of an integral.

		\begin{lemma}
			Suppose that
			\begin{itemize}[label=\textbullet, leftmargin=*]
				\item $\alpha_1\in\ClassK$, and 
				\item $\displaystyle{\int_0^\infty \alpha_1(\norm{x(t)})\drv t < \infty}$ a.s.
			\end{itemize}
			Then $\displaystyle{\Prob{\lim_{t\ra\infty} x(t) = 0} = 1}$.%
			\label{l:asconv}
		\end{lemma}

		\begin{proof}
			Suppose that the claim is false. Then with positive probability there exists some $\eps' > 0$ and a monotone increasing divergent sequence $(s_i)_{i\in\N}$ in $\posR$ such that $\alpha_1(\norm{x(s_i)}) > \eps'$ for all $i$. (Of course, the sequence depends on the event in the aforesaid set of positive probability.) Since $\{f_p\}_{p\in\PSet}$ is a locally Lipschitz and finite family of vector fields, there exists some $\eps'' > 0$ and $L_{\eps''} > 0$ such that $\sup_{\substack{p\in\PSet,\\\norm{x}\in[0, \eps''[}}\norm{f_p(x)} \le L_{\eps''}\norm{x}$. Let $\eps := \eps'\wedge\eps''$. Note that $\fa x\in\R^n\setmin\{0\}$
			\begin{equation}
				\indic{]0, \eps[}(x(t))\left\lvert\frac{\drv\norm{x(t)}}{\drv t}\right\rvert \le L_\eps\norm{x(t)}.
				\label{e:randstab:12}
			\end{equation}
			By the finiteness condition on the integral in the hypothesis, almost surely there exists $T(\eps) > 0$ such that 
			\begin{equation}
				\int_{T(\eps)}^\infty \alpha_1(\norm{x(t)})\drv t < \frac{1}{2}\int_0^{\frac{\ln 2}{L_\eps}}\alpha_1\biggl(\frac{\eps}{2}\epower{-L_\eps s}\biggr)\drv s,
				\label{e:randstab:13}
			\end{equation}
			where the right hand side is a strictly positive quantity since $\alpha_1\in\ClassKinfty$. By construction, $(s_i)_{i\in\N}$ is a monotone increasing divergent sequence with $y(s_i) > \eps$ with positive probability, and therefore there exists $i(\eps)\in\N$ such that $s_{i(\eps)} > T(\eps)$ with strictly positive probability. By continuity of $\norm{\cdot}$ and $x(\cdot)$, there exists an instant $t' > s_{i(\eps)}$ such that $\norm{x(t')} = \eps/2$, also with positive probability. But by~\eqref{e:randstab:12} we have $\norm{x(t)}\in\,]0,\eps[$ for all $t\in\,]t', t'+\frac{\ln 2}{L_\eps}[$, and therefore
			\[
				\int_{t'}^{t'+\frac{\ln 2}{L_\eps}} \alpha_1(\norm{x(t)})\drv t \ge \int_{t'}^{t'+\frac{\ln 2}{L_\eps}}\alpha_1\biggl(\frac{\eps}{2}\epower{-L_\eps(t-t')}\biggr)\drv t
			\]
			with positive probability, which is a contradiction in view of~\eqref{e:randstab:13}. The thesis follows.
		\end{proof}

		\begin{lemma}
			The system~\eqref{e:ssysdef} has the following property: for every $\eps > 0$ there exists $L_\eps > 0$ such that $\norm{x(t)} \le \norm{\xz}\epower{L_\eps t}\quad\fa t\ge 0$ as long as $\norm{x(t)} < \eps$.
			\label{l:detest}
		\end{lemma}

		This is a standard calculation that employs the locally Lipschitz condition on the set of vector fields $\{f_p\}_{p\in\PSet}$.

		\begin{proof}[Proof of Theorem~\ref{t:randstab}]
			To prove that~\eqref{e:ssysdef} is \gasas{}, we need to verify the (SP1) and (SP2) properties in Definition~\ref{d:gasas}.

			From Lemma~\ref{l:Vsigmaest} and Lemma~\ref{l:mgenfun} it follows that $\int_0^\infty \Expec{V_{\sigma(t)}(x(t))}\drv t < \infty$, and by Tonelli's theorem we have
			\[
				\int_0^\infty \Expec{V_{\sigma(t)}(x(t))}\drv t = \Expec{\int_0^\infty V_{\sigma(t)}(x(t))\drv t} < \infty.
			\]
			By hypothesis~\ref{cond:randstab:1} we get $\int_0^\infty \alpha_1(\norm{x(t)})\drv t < \infty$ a.s., and Lemma~\ref{l:asconv} shows that $\lim_{t\ra\infty}\norm{x(t)} = 0$ a.s., which proves (SP2)

			Now we verify (SP1). Fix $\eps > 0$. We know from the (SP2) property proved above that almost surely there exists $T(1, \eps) > 0$ such that $\norm{\xz} < 1$ implies $\sup_{t\ge T(1,\eps)}\norm{x(t)} < \eps$. Select $\delta(\eps) = \min\left\{\eps \epower{-L_\eps T(1, \eps)}, 1\right\}$. By Lemma~\ref{l:detest}, $\norm{\xz} < \delta(\eps)$ implies
			\[
				\norm{x(t)} \le \norm{\xz}\epower{L_\eps t} < \delta(\eps)\epower{L_\eps T(1, \eps)} < \eps \quad \fa t\in[0, T(1, \eps)].
			\]
			Further, the (SP2) property guarantees that with the above choice of $\delta$ and $\xz$, we have $\sup_{t\ge T(1, \eps)} \norm{x(t)} < \eps$ on a set of full measure. Thus, $\norm{\xz} < \delta(\eps)$ implies $\sup_{t \ge 0}\norm{x(t)} < \eps$ a.s. Since $\eps$ is arbitrary, the (SP1) property of~\eqref{e:ssysdef} follows.

			We conclude that~\eqref{e:ssysdef} is \gasas{}
		\end{proof}

		\begin{remark}
			Besides \gasas{}, global asymptotic stability in the mean is another important stability concept. The system~\eqref{e:ssysdef} is said to be \textit{globally asymptotically stable in the mean} (\gasm{}) if the following two properties are simultaneously verified:
			\begin{enumerate}[label={\rm (SM\arabic*)}, align=right, leftmargin=*]
				\item $\fa\eps > 0\quad \therex\wt\delta > 0$ $\norm{\xz} < \wt\delta$ implies $\sup_{t\ge 0}\Expec{\norm{x(t)}} < \eps$;\medskip
				\item $\fa r, \eps' > 0\quad \therex \wt T \ge 0$ $\norm{\xz} < r$ implies $\sup_{t\ge \wt T}\Expec{\norm{x(t)}} < \eps'$.
			\end{enumerate}
			We have seen that under the hypotheses of Theorem~\ref{t:randstab} we have convergence of $\bigl(\alpha_1(\norm{x(t)})\bigr)_{t\ge 0}$ to $0$ in $L_1$. If $\alpha_1$ is convex, then this immediately gives $\lim_{t\ra\infty}\Expec{\norm{x(t)}} = 0$ via Jensen's inequality,\footnote{We recall Jensen's inequality~\cite[p.~80]{ref:billingsleybk}: if $X$ is an integrable real-valued random variable on $(\Omega, \sigalg, \mathsf P)$, and if $\phi:\R\lra\R$ is a convex function, then $\phi\big(\Expec{\left|X\right|}\big) \le \Expec{\vphantom{\big(}\phi(\left|X\right|)}$.} which is (SM2). Lemma~\ref{l:Vsigmaest} and Lemma~\ref{l:mgenfun} show that $\sup_{t\ge 0}\Expec{\alpha_1(\norm{x(t)})} \le \alpha_2(\norm{\xz})(S+1)$, which under the assumption that $\alpha_1$ is convex, implies (SM1). In practice the Lyapunov functions are usually taken to be polynomial powers of the state $x$ and $\alpha_1$ is convex. Otherwise a property analogous to \gasm{} still holds, with $\alpha_1(\norm{x(t)})$ replacing $\norm{x(t)}$ everywhere, and therefore depends on the choice of $\alpha_1$.\RemarkEnd
		\end{remark}

		\begin{remark}
			With a slight modification in the hypotheses of Theorem~\ref{t:randstab} we can employ standard results in martingale theory to conclude almost sure global asymptotic convergence of $x(\cdot)$. Indeed, if the condition~\eqref{cond:randstab:4} is strengthened to the conditional version
			\[
			\Prob{N_\sigma(t+s) - N_\sigma(t) = k\mid\sigalg_t} \le \epower{-\wt\lambda(s)}\frac{\bigl(\ol\lambda s\bigr)^k}{k!}
			\]
			for $k\in\N\cup\{0\}, \;t\ge 0,\; s > 0$, then a calculation in the spirit of Lemma~\ref{l:Vsigmaest} shows that $\bigl(V_{\sigma(t)}(x(t))\bigr)_{t\ge 0}$ is a supermartingale. Lemma~\ref{l:Vsigmaest} also shows that $\lim_{t\ra\infty}\Expec{V_{\sigma(t)}(x(t))} = 0$, which implies that the aforesaid process is a potential. A standard result in martingale theory (e.g.,~\cite[p.~18, Problem~3.16]{ref:karatshreveStochCalc}) now implies that the process $\bigl(V_{\sigma(t)}(x(t))\bigr)_{t\ge 0}$ converges to $0$ a.s. Considering hypothesis~\ref{cond:randstab:1} of Theorem~\ref{t:randstab} we conclude that $\bigl(\norm{x(t)}\bigr)_{t\ge 0}$ converges to $0$ a.s.\RemarkEnd
		\end{remark}

		\subsection{Markov jump systems}
		We note that hypothesis~\ref{cond:randstab:4} of Theorem~\ref{t:randstab} stipulates that $\fa t\in\posR$ the tail of the probability mass function of the random variable $N_\sigma(t)$ is majorized (i.e., stochastically dominated) by the probability mass function of a ``maximally'' switching jump-stochastic process. This hypothesis can be verified, in particular, if $\sigma$ is the state of a continuous-time Markov chain, with a given generator matrix $Q = [q_{ij}]_{\cardP\times\cardP}^{\vphantom{T}}$ and a given initial probability distribution $\pi^\circ$ (recall that $\cardP$ is the number of elements of $\PSet$); we denote this by $\sigma\sim(\pi^\circ, Q)$. Lemma~\ref{l:mcest} and Corollary~\ref{c:randstabmc} make this statement precise.

		Let us recall some basic facts about continuous-time Markov chains, see e.g.~\cite{ref:norrisbk} for further details. Associated with the Markov chain $\sigma\sim(\pi^\circ, Q)$ is the Kolmogorov forward equation
		\[
			\dot P(t) = P(t) Q,\qquad P(0) = \idmat{\cardP}, \quad t\ge 0,
		\]
		where $\idmat{\cardP}$ is the $\cardP$-dimensional identity matrix; the probability (row) vector at any time $t\ge 0$ is given by $\pi(t) = \pi^\circ P(t)$.
		We need the following two facts:
		\begin{enumerate}[label={\rm (MC\arabic*)}, leftmargin=*, align=right]
			\item The generator matrix $Q = [q_{i j}]_{\cardP\times\cardP}$ satisfies $\big(q_{ij} \ge 0\quad i\neq j\big)$, and $\bigl(\sum_{j\in\PSet\smallsetminus\{i\}} q_{ij} = -q_{ii}\bigr)$ for $i, j\in\PSet$.
			\item $\Prob{\sigma(t+h) = j\,\left|\vphantom{\big(}\right.\,\sigma(t) = i} = \delta_{ij} + q_{ij}h + o(h)$ for $h > 0$, and $\delta_{ij}$ is the Kronecker delta. This is known as the infinitesimal description of a continuous-time Markov chain.
		\end{enumerate}

		We define
		\begin{equation}
			\ol q := \max\big\{\lvert q_{ii}\rvert \,\big|\,i\in\PSet\big\},\quad \wt q := \max \big\{q_{ij}\,\big|\,i, j\in\PSet\big\}.
			\label{e:Qparamdef}
		\end{equation}

		The following lemma establishes a coarse estimate of the probability mass function of $N_\sigma(t)$; a finer estimate can be obtained from the Levy formula, but for our purposes the one derived below is sufficient.

		\begin{lemma}
			Suppose that $\sigma\sim(\pi^\circ, Q)$ is a Markov chain. Then $\fa t\in\posR$, we have $\Prob{N_\sigma(t) = k} \le \epower{-\wt qt}(\ol qt)^k/k!\quad \fa k\in\N\cup\{0\}$.
			\label{l:mcest}
		\end{lemma}
		\begin{proof}
			For $t\in\posR$ and $k\in\N\cup\{0\}$, define $\eta_k(t) := \Prob{N_\sigma(t) = k}$. For $h > 0$ sufficiently small, $\fa k\in\N\cup\{0\}$,
			\begin{equation}
				\eta_k(t+h) = \sum_{i = 0}^k \Prob{\!N_\sigma(t+h)-N_\sigma(t) = i\!}\Prob{\!N_\sigma(t) = k-i\!}.
				\label{e:etatph}
			\end{equation}
			By the infinitesimal description of a Markov chain (MC2),
			\begin{align}
				& \Prob{N_\sigma(t+h) - N_\sigma(t) = 0}
				 \le 1-\wt qh + o(h),
				\label{e:diffiszero}
			\end{align}
			and
			\begin{align}
				& \Prob{N_\sigma(t+h)-N_\sigma(t) = 1}
				 \le \ol qh + o(h).
				\label{e:diffisone}
			\end{align}
			For all natural numbers $k\ge 2$, (MC2) shows that
			\begin{equation}
				\Prob{N_\sigma(t+h)-N_\sigma(t) = k} = o(h).
				\label{e:diffmorethantwo}
			\end{equation}
			Using~\eqref{e:diffiszero},~\eqref{e:diffisone} and~\eqref{e:diffmorethantwo}, we continue the calculation in~\eqnref{e:etatph}:
			\[
				\eta_k(t+h) \le (1-\wt qh + o(h))\eta_k(t) + (\ol qh + o(h))\eta_{k-1}(t) + o(h),
			\]
			which leads to
			\[
				\frac{\eta_k(t+h) - \eta_k(t)}{h} \le -\wt q\eta_k(t) + \ol q\eta_{k-1}(t) + O(h).
			\]
			Taking limits with $h\downarrow 0$, the following differential inequality is obtained:
			\[
				\dot\eta_k(t) \le -\wt q\eta_k(t) + \ol q\eta_{k-1}(t),\qquad \eta_k(0) = 0,\quad\fa k\in\N.
			\]
			(We have identical differential inequalities starting with $t > 0$ and $h < 0$ sufficiently small.) A similar analysis yields
			\[
				\dot \eta_0(t) \le -\wt q\eta_0(t),\qquad \eta_0(0) = 1.
			\]
			In matrix notation, the set of differential inequalities involving $\dot\eta_k,\; k\in\N\cup\{0\}$, stands as:
			\begin{equation}
				\begin{bmatrix}
					\dot \eta_0\\
					\dot \eta_1\\
					\dot \eta_2\\
					\vdots
				\end{bmatrix}
				\le
				\begin{bmatrix}
					-\wt q & 0 & 0 & \cdots\\
					\ol q & -\wt q & 0 & \cdots\\
					0 & \ol q & -\wt q & \cdots\\
					\vdots & \vdots & \vdots & \ddots
				\end{bmatrix}
				\begin{bmatrix}
					\eta_0\\
					\eta_1\\
					\eta_2\\
					\vdots
				\end{bmatrix},\quad
				\begin{bmatrix}
					\eta_0(0)\\
					\eta_1(0)\\
					\eta_2(0)\\
					\vdots
				\end{bmatrix}
				=
				\begin{bmatrix}
					1\\
					0\\
					0\\
					\vdots
				\end{bmatrix},
				\label{e:mainineq}
			\end{equation}
			where the ``$\le$'' is interpreted component-wise. Clearly, $\eta_0(t) \le \epower{-\wt q t}, \; t\ge 0,$ satisfies the first differential inequality. We claim that 
			\begin{equation}
				\eta_k(t) \le \epower{-\wt q t}(\ol qt)^k/k!\qquad \fa t\ge 0\quad \fa k\in\N
				\label{e:etakt}
			\end{equation}
			is a solution to~\eqref{e:mainineq}. Indeed, for $k = 1$ we have $\dot \eta_1 \le \ol q\eta_0 - \wt q\eta_1 \le \ol q\epower{-\wt qt} - \wt q\eta_1$, which leads to 
			\[
				\epower{\wt qt}\eta_1(t) \le \epower{\wt qt}\eta_1(0) + \ol q\int_{0}^t\drv s;
			\]
			hence $\eta_1\epower{\wt qt} \le (\ol qt)$ (in view of $\eta_1(0) = 0$), yielding $\eta_1(t) \le (\ol qt)\epower{-\wt qt},\; t\ge 0$. Having verified the claim for $k = 1$,
			an induction argument shows that for arbitrary $j$
			\[
				\epower{\wt qt}\eta_{j+1}(t) \le \epower{\wt qt}\eta_{j+1}(0) + \ol q\int_{0}^t\frac{(\ol qs)^j}{j!}\drv s;
			\]
			hence $\eta_{j+1}\epower{\wt qt} \le (\ol qt)^{j+1}/{(j+1)}!$ (in view of $\eta_{j+1}(0) = 0$), yielding $\eta_{j+1}(t) \le \epower{-\wt qt}(\ol qt)^{j+1}/{(j+1)}!,\; t\ge 0$.
			In view of the definition of $\eta_k(t)$, the thesis of the Lemma follows.
		\end{proof}

		\begin{corollary}
			Consider the system~\eqref{e:ssysdef}, and let $\ol q, \wt q$ be defined by~\eqref{e:Qparamdef}. Suppose that $\sigma\sim(\pi^\circ, Q)$ is a Markov chain, and that there exist continuously differentiable functions $V_p:\R^n\lra\posR,\;p\in\PSet$, functions $\alpha_1, \alpha_2\in\ClassKinfty$, and a real number $\mu > 1$, such that
			\begin{enumerate}[label={\rm (\roman*)}, align=right, leftmargin=*]
				\item $\alpha_1(\norm x) \le V_p(x) \le \alpha_2(\norm x)\quad \fa x\in\R^n,\;\fa p\in\PSet$;\label{cond:randstabmc:1}\medskip
				\item $\displaystyle{\frac{\partial V_p}{\partial x}f_p(x) \le -\lambda_\circ V_p(x)}\quad\fa x\in\R^n,\; \fa p\in\PSet$;\label{cond:randstabmc:2}\medskip
				\item $V_{p_1}(x) \le \mu V_{p_2}(x)\quad\fa x\in\R^n,\; \fa p_1, p_2\in\PSet$;\label{cond:randstabmc:3}\medskip
				\item $\mu < (\lambda_\circ+\wt q)/\ol q$.\label{cond:randstabmc:4}
			\end{enumerate}
			Then~\eqref{e:ssysdef} is \gasas{}
			\label{c:randstabmc}
		\end{corollary}
		\begin{proof}
			Follows directly from Lemma~\ref{l:mcest}, and Theorem~\ref{t:randstab} with $M = 0$ in hypothesis~\ref{cond:randstab:4}.
		\end{proof}

	\section{Stabilization of Randomly Switched Control Systems}
		In this section we establish a method for designing controllers that ensure almost sure global asymptotic stability of control-affine randomly switched systems in closed loop. The method may be viewed as an application of our results in~\S\ref{s:mainres}. We assume that at each instant of time $t$ the state $\sigma(t)\in\PSet$ of the random switching signal is perfectly known to the controller. 

		Consider the affine in control switched system:
        \begin{equation}
            \dot x = f_\sigma(x) + \sum_{i = 1}^mg_{\sigma, i}(x) u_i, \qquad x(0) = \xz, \quad t\ge 0,
            \label{e:ssysdefcon}
        \end{equation}
		where $x\in\R^n$ is the state, $u_i,\;i = 1, \ldots, m$ are the control inputs, $u_i\in\R$, $f_p$ and $g_{p, i}$ are smooth vector fields on $\R^n$, with $f_p(0) = 0, g_{p, i}(0) = 0$, for each $p\in\PSet, i\in\{1, \ldots, m\}$. With a feedback control function $\ol u_\sigma(x) = \left[u_{\sigma, 1}(x), \ldots, u_{\sigma, m}(x)\right]\transp$, the closed loop system stands as:
		\begin{equation}
			\dot x = f_\sigma(x) + \sum_{i = 1}^m g_{\sigma, i}(x) \ol u_{\sigma, i}(x), \quad x(0) = \xz, \;\; t\ge 0.
			\label{e:ssysdefconcl}
		\end{equation}
		Our objective is to select the control function $\ol u_\sigma$ so that~\eqref{e:ssysdefconcl} is \gasas{} Let the switching signal $\sigma$ be a stochastic process as defined in \S\ref{s:Prelims}, and let $\xz\neq 0$.

		A universal formula for stabilization of control-affine nonlinear systems was first constructed in~\cite{ref:sontagunivformula}, for the control taking values in $\R^m$.

		\begin{theorem}
			Consider the system~\eqref{e:ssysdefcon}. Suppose that $\sigma$ satisfies hypothesis~\ref{cond:randstab:4} of {\rm Theorem~\ref{t:randstab}}, and there exists a family of continuously differentiable functions $V_p:\R^n\lra\posR,\;{p\in\PSet}$, such that
			\begin{enumerate}[label={\rm (C\arabic*)}, align=right, leftmargin=*]
				\item hypothesis~\ref{cond:randstab:1} of {\rm Theorem~\ref{t:randstab}} holds;\label{cond:u:1}\medskip
				\item hypothesis~\ref{cond:randstab:3} of {\rm Theorem~\ref{t:randstab}} holds;\label{cond:u:4}\medskip
				\item $\therex\lambda_\circ > 0$ such that $\fa x\in\R^n\setmin\{0\}$ and $\fa p\in\PSet$\label{cond:u:2}
				\[
					\inf_{u\in\R^m}\left\{\!\LieD{f_p}{V_p(x)}\! + \!\lambda_\circ V_p(x) \!+ \!\sum_{i = 1}^m \!u_i\LieD{g_{p, i}}{V_p(x)}\right\} < 0;
				\]
				\item $\fa \eps > 0 \;\therex \delta > 0$ such that if $x(\neq 0)$ satisfies $\norm{x} < \delta$, then $\therex u\in\R^m,\; \norm u < \eps$, such that $\fa p\in\PSet$\label{cond:u:3}
				\[
					\LieD{f_p}{V_p} + \sum_{i = 1}^m u_i\cdot\LieD{g_{p, i}}{V_p} \le -\lambda_\circ V_p;
				\]
				\item hypothesis~\ref{cond:randstab:5} of {\rm Theorem~\ref{t:randstab}} holds.
			\end{enumerate}
			Then the feedback control
			\[
				\ol u_\sigma(x) = [k_{\sigma, 1}(x), \ldots, k_{\sigma, m}(x)]\transp,
			\]
			where
			\begin{subequations}
			\begin{align}
				k_{p, i}(x) & := \displaystyle{-\LieD{g_{p, i}}{V_p}(x)\cdot\varphi\!\left(\ol W_p(x), \wt W_p(x)\right)} \\
				\ol W_p(x) & := \LieD{f_p}{V_p}(x) + \lambda_\circ V_p(x),\\
				\wt W_p(x) & := \sum_{i = 1}^m\left(\LieD{g_{p, i}}{V_p}(x)\right)^2,\\
			\intertext{and}
				\varphi(a, b) & := 
				\begin{cases}
					\displaystyle{\frac{a + \sqrt{a^2 + b^2}}{b}}\quad & \text{if }b\neq 0,\\
					0 & \text{otherwise,}
				\end{cases}
				\label{e:phidef}
			\end{align}
			\label{e:udefs}
			\end{subequations}
			renders~\eqref{e:ssysdefconcl} \gasas{}
			\label{t:u}
		\end{theorem}
		\begin{proof}
			The proof relies on the construction of the universal formula in~\cite{ref:sontagunivformula}.
			Fix $t\in\posR$. If $x\neq 0$, applying the definition of $\varphi$, we get
			\begin{align*}
				& \LieD{f_{\sigma(t)}}{V_{\sigma(t)}}(x) + \sum_{i = 1}^m k_{\sigma(t), i}(x)\LieD{g_{\sigma(t), i}}{V_{\sigma(t)}}(x)\\
				& = \LieD{f_{\sigma(t)}}{V_{\sigma(t)}}(x) - \wt W_{\sigma(t)}(x)\!\cdot\!\varphi\!\left(\!\ol W_{\sigma(t)}(x),\!\left(\wt W_{\sigma(t)}(x)\right)^{\!2}\right)\\
				& = -\lambda_{\circ} V_{\sigma(t)}(x) -\sqrt{\left(\LieD{f_{\sigma(t)}}{V_{\sigma(t)}}(x)\right)^2 + \left(\wt W_{\sigma(t)}(x)\right)^2}\\
				& < -\lambda_{\circ} V_{\sigma(t)}(x).
			\end{align*}
			Since $t$ is arbitrary, we conclude that the above inequality holds for all $t\in\posR$. Note that by (C3), if for any $p\in\PSet$, $x\in\bigcap_{i = 1}^m \ker\left(\LieD{g_{p, i}}{V_p}\right)$, we automatically have $\LieD{f_{\sigma(t)}}{V_{\sigma(t)}}(x) + \lambda_{\circ} V_{\sigma(t)}(x) < 0$. (C4) is the small control property, ensuring continuity of the control function at $0$ for each fixed index $p$; this guarantees the existence of a unique local solution to the switched system.

			The above arguments, in conjunction with (C1) and (C2) enable us to conclude that the family $(V_p)_{p\in\PSet}$ satisfies hypotheses~\ref{cond:randstab:1},~\ref{cond:randstab:2} and~\ref{cond:randstab:3} of Theorem~\ref{t:randstab} for the closed loop system~\eqref{e:ssysdefconcl}. (C5) ensures that hypothesis~\ref{cond:randstab:5} of Theorem~\ref{t:randstab} holds for~\eqref{e:ssysdefconcl}. Since $\sigma$ satisfies hypothesis~\ref{cond:randstab:4} of Theorem~\ref{t:randstab}, it follows from Theorem~\ref{t:randstab} applied to~\eqref{e:ssysdefconcl}, that~\eqref{e:ssysdefconcl} is \gasas{}
		\end{proof}

	\section{Conclusion and Further Work}
		We have provided sufficient conditions for almost sure stability of randomly switched systems, together with control strategies for almost sure stabilization for systems with control inputs. It may be possible to improve upon the proposed results by utilizing the jump destinations of the switching signal, and in the case of Markov chains, its graph and the associated transition probability matrix. Stabilization of randomly switched systems with control inputs without perfect knowledge of $\sigma$ is a nontrivial and important issue. Input-to-state stability properties, existence and uniqueness of invariant measures, and other asymptotic properties of randomly switched systems are interesting avenues for future research. Results on these will be reported elsewhere.

	\section*{Acknowledgments}
		This study was motivated by a conversation with P.~K.~Jha and L.~Vu.

\bigskip

	\bibliographystyle{siam}
	\bibliography{../references}

\begin{thebibliography}{10}

\bibitem{ref:basarMinMaxSw}
{\sc T.~Ba{\c s}ar}, {\em Minimax control of switching systems under sampling},
  Systems \& Control Letters, 25 (1995), pp.~315--325.

\bibitem{ref:battilottiDwellTimeMJS}
{\sc S.~Battilotti and A.~D. Santis}, {\em Dwell time controllers for
  stochastic systems with switching {M}arkov chain}, Automatica, 41 (2005),
  pp.~923--934.

\bibitem{ref:billingsleybk}
{\sc P.~Billingsley}, {\em Probability and Measure}, Wiley-Interscience, 3~ed.,
  1995.

\bibitem{ref:bolzernCDC04}
{\sc P.~Bolzern, P.~Colaneri, and G.~D. Nicolao}, {\em On almost sure stability
  of discrete-time {M}arkov jump linear systems}, in Proceedings of the 43rd
  Conference on Decision and Control, 2004, pp.~3204--3208.

\bibitem{ref:lygerosGenStochHybSys}
{\sc M.~L. Bujorianu and J.~Lygeros}, {\em General stochastic hybrid systems:
  modeling and optimal control}, in Proceedings of the 43rd IEEE Conference on
  Decision and Control, 2004, pp.~1872--1877.

\bibitem{ref:davisMarkovModelsOptbk}
{\sc M.~H.~A. Davis}, {\em Markov Models and Optimization}, Chapman \& Hall,
  London, 1993.

\bibitem{ref:fengMarkovJump92}
{\sc X.~Feng, K.~A. Loparo, Y.~Ji, and H.~J. Chizeck}, {\em Stochastic
  stability properties of jump linear systems}, IEEE Transactions on Automatic
  Control, 37 (1992), pp.~38--53.

\bibitem{ref:filippovbk}
{\sc A.~F. Filippov}, {\em Differential Equations with Discontinuous Righthand
  Sides}, vol.~18 of Mathematics and Its Applications, Kluwer Academic
  Publishers, Dordrecht, 1988.

\bibitem{ref:khasminskiiStochastic}
{\sc R.~Z. Ha{\'s}minskii}, {\em Stochastic Stability of Differential
  Equations}, Sijthoff \& Noordhoff, Alphen aan den Rijn - Germantown, 1980.

\bibitem{ref:hespanhaCommNet}
{\sc J.~P. Hespanha}, {\em A model for stochastic hybrid systems with
  application to communication networks}, Nonlinear Analysis, 62 (2005),
  pp.~1353--1383.

\bibitem{ref:hespanhaAvgDwellTime}
{\sc J.~P. Hespanha and A.~S. Morse}, {\em Stability of switched systems with
  average dwell-time}, in Proceedings of the 38th IEEE Conference on Decision
  and Control, vol.~3, 1999, pp.~2655--2660.

\bibitem{ref:meynODESpecMarkovOp}
{\sc J.~Huang, I.~Kontoyiannis, and S.~P. Meyn}, {\em The {ODE} method and
  spectral theory of {M}arkov operators}, in Stochastic Theory and Control
  (Lawrence, KS, 2001), vol.~280 of Lecture Notes in Control and Information
  Sciences, Springer, Berlin, 2002, pp.~205--221.

\bibitem{ref:jichizeck90}
{\sc Y.~Ji and H.~J. Chizeck}, {\em Controllability, stabilizability, and
  continuous-time {M}arkovian jump linear quadratic control}, IEEE Transactions
  on Automatic Control, 35 (1990), pp.~777--788.

\bibitem{ref:karatshreveStochCalc}
{\sc I.~Karatzas and S.~E. Shreve}, {\em Brownian {M}otion and {S}tochastic
  {C}alculus}, vol.~113 of Graduate Texts in Mathematics, Springer-Verlag, New
  York, 2~ed., 1991.

\bibitem{ref:kozinsurvey}
{\sc F.~Kozin}, {\em A survey of stability of stochastic systems}, Automatica,
  5 (1969), pp.~95--112.

\bibitem{ref:lakshmikanthamDifferentialIntegralInequalities1}
{\sc V.~Lakshmikantham and S.~Leela}, {\em Differential and Integral
  Inequalities: Theory and Application}, vol.~1, Academic Press, New York -
  London, 1969.

\bibitem{ref:liberzonbk}
{\sc D.~Liberzon}, {\em Switching in Systems and Control}, Systems \& Control:
  Foundations \& Applications, Birkh{\"a}user, Boston, 2003.

\bibitem{ref:maoexpstabstocdelay}
{\sc X.~Mao}, {\em Exponential stability of stochastic delay interval systems
  with {M}arkovian switching}, IEEE Transactions on Automatic Control, 47
  (2002), pp.~1604--1612.

\bibitem{ref:norrisbk}
{\sc J.~R. Norris}, {\em Markov Chains}, Cambridge University Press, Cambridge,
  1997.

\bibitem{ref:sontagunivformula}
{\sc E.~D. Sontag}, {\em A universal construction of {A}rtstei{n's} theorem on
  nonlinear stabilization}, Systems \& Control Letters, 13 (1989),
  pp.~117--123.

\bibitem{ref:xiong05}
{\sc J.~Xiong, J.~Lam, H.~Gao, and D.~W.~C. Ho}, {\em On robust stabilization
  of {M}arkovian jump systems with uncertain switching probabilities},
  Automatica, 41 (2005), pp.~897--903.

\end{thebibliography}

\end{document}